\documentclass[12pt]{amsart}
\usepackage{amsmath}
\usepackage{amsxtra}
\usepackage{amscd}
\usepackage{amsthm}
\usepackage{amsfonts}
\usepackage{amssymb}
\usepackage{eucal}

\numberwithin{equation}{subsection}

\begin{document}

\newtheorem{thm}[subsection]{Theorem}
\newtheorem{prop}[subsection]{Proposition}
\newtheorem{lem}[subsection]{Lemma}
\newtheorem{cor}[subsection]{Corollary}
\newtheorem{rmk}[subsection]{Remark}
\newtheorem{conv}[subsection]{Convention}
\newtheorem{conj}[subsection]{Conjecture}
\newtheorem{defn}[subsection]{Definition}
\newtheorem{ass}[subsection]{Assumption}
\newtheorem{notation}[subsection]{Notation}
\newtheorem{example}[subsection]{Example}
\newtheorem{caution}[subsection]{Caution}
\newcommand{\thmref}[1]{Theorem~\ref{#1}}
\newcommand{\defref}[1]{Definition~\ref{#1}}
\newcommand{\lemref}[1]{Lemma~\ref{#1}}
\newcommand{\propref}[1]{Proposition~\ref{#1}}
\newcommand{\corref}[1]{Corollary~\ref{#1}}
\newcommand{\rmkref}[1]{Remark~\ref{#1}}
\newcommand{\assref}[1]{Assumption~\ref{#1}}

\newcommand{\nc}{\newcommand}
\nc{\ol}{\overline}
\nc{\on}{\operatorname}

\nc{\BA}{\mathbb A}
\nc{\BC}{\mathbb C}
\nc{\BD}{\mathbb D}
\nc{\BG}{\mathbb G}
\nc{\BH}{\mathbb H}
\nc{\BR}{\mathbb R}
\nc{\BZ}{\mathbb Z}

\nc{\CB}{\mathcal B}
\nc{\CC}{\mathcal C}
\nc{\CF}{\mathcal F}
\nc{\CJ}{\mathcal J}
\nc{\CK}{\mathcal K}
\nc{\CL}{\mathcal L}
\nc{\CM}{\mathcal M}
\nc{\CN}{\mathcal N}
\nc{\CO}{\mathcal O}
\nc{\CP}{\mathcal P}
\nc{\CR}{\mathcal R}
\nc{\CS}{\mathcal S}
\nc{\CT}{\mathcal T}
\nc{\CU}{\mathcal U}
\nc{\CV}{\mathcal V}
\nc{\CW}{\mathcal W}
\nc{\CX}{\mathcal X}

\nc{\frakb}{\mathfrak b}
\nc{\fraki}{\mathfrak i}
\nc{\frakj}{\mathfrak j}
\nc{\frakn}{\mathfrak n}

\nc{\risom}{\buildrel {\sim}\over {\rightarrow}}
\nc{\lisom}{\buildrel {\sim}\over {\leftarrow}}

\nc{\Sym}{{\on{Sym}}}
\nc{\vect}{{\on{Vect}}}
\nc{\Hom}{{\on{Hom}}}
\nc{\Lie}{{\on{Lie}}}
\nc{\Id}{{\on{Id}}}
\nc{\IC}{{\on{IC}}}
\nc{\inthom}{\mathcal Hom}
\nc{\mult}{{\on{mult}}}

\nc{\verdier}{\mathcal D}


\title{Morse theory and tilting sheaves}

\author{David Nadler}

\address{Department of Mathematics\\ 
Northwestern University\\ Evanston, IL 60208}
\email{nadler@math.northwestern.edu}


\maketitle

\begin{center} 
\em for Bob
\end{center}



\section{Introduction}

Tilting perverse sheaves arise in geometric representation theory as natural bases
for many categories.
(See~\cite{BBM04} or Section~\ref{secsh} below for definitions and some background results.)
For example, for category $\CO$ of a complex reductive group $G$, tilting modules correspond 
to tilting sheaves  
on the Schubert stratification $\CS$ of the
flag variety $\CB$.
In this setting,
there are many characterizations of tilting objects involving 
their categorical structure. 
(See~\cite{Ringel91} and~\cite{SoeCharForm98}.)
In this paper, we propose a new geometric
construction of tilting sheaves. Namely, we show that they naturally
arise via the Morse theory
of sheaves on the opposite Schubert stratification $\CS^{opp}$.

In the context of the flag variety $\CB$,
our construction takes the following form.
Let $\CF$ be a complex of sheaves on $\CB$ with bounded constructible cohomology. 
We will refer to such an object as simply a sheaf.

Given a point $p\in\CB$, and the germ $F$ at $p$ of a complex analytic function, 
we have the corresponding local Morse group (or vanishing cycles) 
$M^*_{p,F}(\CF)$ which measures how
$\CF$ changes at $p$ as we move in the family defined by $F$. 
The local Morse groups play a central role in the theory of perverse sheaves. 
For an arbitrary sheaf on $\CB$, its local Morse groups are graded by cohomological degree.
The perverse sheaves on $\CB$ are characterized by the fact that their
local Morse groups 
(with respect to sufficiently generic functions $F$)
are concentrated in degree $0$. In other words,
the local Morse groups are $t$-exact functors (with respect to the perverse $t$-structure),
and any sheaf on which they are concentrated in degree $0$ is in turn perverse.

If we restrict to sheaves $\CF$ on the Schubert stratification
$\CS^{opp}$,
then for any function $F$ at a point $p$,
it is possible to find a sheaf $\CM_{p,F}$ on a neighborhood of $p$
which represents the functor $M^*_{p,F}$ (though $\CM_{p,F}$ is not constructible
with respect to $\CS^{opp}$.)
To be precise, for any $\CF$ on $\CS^{opp}$, we have a natural isomorphism
$$
M^*_{p,F}(\CF)\simeq R\Hom^*_{D(\CB)}(\CM_{p,F},\CF).
$$
Since this may be thought of as an integral identity, 
we call the sheaf $\CM_{p,F}$ a Morse kernel.

Next we choose a one-parameter subgroup $\BC^\times \subset G$ so that
for the resulting
$\BC^\times$-action on $\CB$, the strata of the Schubert stratification $\CS$
are its ascending manifolds, and the strata of 
the opposite Schubert stratification $\CS^{opp}$ are its descending manifolds.
In other words, the action has finitely many fixed points, and for $z\in\BC^\times$,
the action of $z$ contracts the strata of $\CS$ to the fixed points as $z\to\infty$,
and contracts the strata of $\CS^{opp}$ to the fixed points as $z\to 0$.
Then we may take a Morse kernel $\CM_{p,F}$ and flow it down along
the $\BC^\times$-orbits to obtain a limit sheaf (or nearby cycles)
$$
\CM_{p,F}\stackrel{z\to 0}{\longrightarrow} \Psi(\CM_{p,F}).
$$
The limit is a sheaf on $\CB$ with, {\em a priori}, no particular shape. 
Our main result is the following. 
(See~\cite{ENV04} for another application of such limits.)

\begin{thm}
For a Morse kernel $\CM_{p,F}$, the shifted limit $T_{p,F}=\Psi(\CM_{p,F})[\dim_\BC\CB]$ 
is a tilting perverse sheaf on the Schubert stratification $\CS$. For a sheaf $\CF$ on the opposite Schubert stratification $\CS^{opp}$,
we have a natural isomorphism
$$
R\Hom^*_{D(\CB)}(\CM_{p,F},\CF)\simeq R\Hom^*_{D(\CB)}(\Psi(\CM_{p,F}),\CF).
$$
\end{thm}

It is easy to see that all tilting sheaves on $\CS$ arise as summands of such limits.
We mention a possible source of confusion.
It is well known that such limits (or nearby cycles) take perverse sheaves to perverse sheaves.
But the Morse kernels are by no means perverse (they are not even $\BC$-constructible.)
Nevertheless, the above theorem asserts that by taking their limits along
the $\BC^\times$-orbits we obtain perverse sheaves.
The delicate part of the proof is seeing the constructibility of the limit.
Once we confirm Morse kernels are constructible 
with respect to stratifications transverse to $\CS^{opp}$, 
we use Bekka's version of stratifications~\cite{Bekka91} to prove the following.

\begin{prop}
If a sheaf $\CM$ is constructible with respect to a stratification transverse to $\CS^{opp}$,
then the limit $\Psi(\CM)$ is constructible with respect to $\CS$.
\end{prop}

To connect with traditional notions of representability, 
consider the opposite unipotent subgroups $N,N^{opp}\subset G$
whose orbits in $\CB$ are the strata of $\CS,\CS^{opp}$ respectively.
Since $N,N^{opp}$ are contractible, the $N,N^{opp}$-equivariant derived categories
are the same as the $\CS,\CS^{opp}$-constructible derived categories
respectively. Thus 
we have
the equivariantization functors
$$
\Gamma^+:D(\CB)\to D_{\CS}(\CB)
\qquad
\gamma^-:D(\CB)\to D_{\CS^{opp}}(\CB)
$$
which are right and left adjoint respectively to the forgetful functor. 
If we restrict $\Gamma^+,\gamma^-$
to $\CS^{opp},\CS$-constructible sheaves respectively, then they become inverse
equivalences of categories.

Using our main theorem, we see that for a Morse kernel $\CM_{p,F}$,
the equivariantization $\gamma^-$ commutes with the limit $\Psi$.
It follows that 
$\gamma^-(\CM_{p,F})$ is a projective perverse sheaf 
on $\CS^{opp}$. 
Thus we have
a triangle of objects which represent the local Morse group 
$$
\begin{array}{ccccc}
&& \CM_{p,F} && \\
& \Psi \swarrow & & \searrow  \gamma^- & \\
T_{p,F}[-\dim_\BC \CB] & & \stackrel{\gamma^-}{\longrightarrow} & & P^{opp}_{p,F}
\end{array}
$$
At the top, we have the Morse kernel $\CM_{p,F}$ which by construction
represents $M_{p,F}$. To the left, we have the shifted tilting
sheaf $T_{p,F}[-\dim_\BC \CB]$ on $\CS^{}$ to which $\CM_{p,F}$ flows. By the above theorem, 
$T_{p,F}[-\dim_\BC \CB]$
also represents $M_{p,F}$. Finally, to the right, we have the projective
perverse sheaf $P^{opp}_{p,F}$ 
on $\CS^{opp}$ which represents $M_{p,F}$ both in the abelian and derived categories.

An immediate corollary is the following multiplicity formula for the limit tilting sheaf
$T_{p,F}$
in terms of the minimal tilting sheaves:
$$
\mult_{T_w}(T_{p,F})=\dim M_{p,F}(\IC^{opp}_w).
$$
Here $T_w$ is the minimal tilting sheaf on the 
Schubert stratum corresponding to $w$ in the Weyl group of $G$,
and $\IC^{opp}_w$ is the intersection cohomology sheaf on the opposite stratum.

As another consequence, we see that
the equivariantization functors exchange the shifted minimal tilting and projective perverse
sheaves on the opposite stratifications
$$
\gamma^-(T_w[-\dim_\BC \CB])\simeq P^{opp}_w
\qquad
\Gamma^+(P^{opp}_w)\simeq T_w[-\dim_\BC \CB]
$$
This statement should be compared to the fact (see \cite{BBM04}) 
that for perverse sheaves on $\CS$ alone, 
the ``longest interwiner" exhanges the minimal tilting and projective perverse sheaves.
The above statement is closely related to this, but continues to make sense in other contexts
such as the affine setting. 


\subsection{Organization of paper}

In what follows, we work in a more general context than discussed in the introduction above.
Namely, we do not restrict ourselves to flag varieties, but rather consider any space with
the necessary properties for our main construction to be easily understood.
This is explained in Section~\ref{secstrat} along with other background material.
In Section~\ref{secsh}, we briefly recall the definition of perverse sheaves and tilting sheaves
and their basic properties.
In Section~\ref{secmorse}, we discuss Morse groups and Morse kernels on stratified spaces.
Section~\ref{secmain} contains our main result and its proof.
Finally, in Section~\ref{secflag}, we return to the case of flag varieties.


\subsection{Acknowledgments}
This paper was inspired by an ongoing project with David Ben-Zvi
concerning real groups and the geometric Langlands program.
It is a pleasure to thank him for his interest and many enlightening conversations. 

I would also like to thank Mark Goresky, Tom Nevins, and Kari Vilonen for useful comments.

Finally, I would like to comment on the debt this work owes to Bob MacPherson.
Through his own work,
he has played a fundamental role in developing 
the main objects appearing here, including
stratified spaces, $\BC^\times$-actions, Morse theory, perverse sheaves, and nearby cycles.
His contributions to these subjects continue to open new avenues of exploration.
He has also had a tremendous influence 
through his generous support and enthusiasm for the work of others.
It is a pleasure to devote this paper to Bob as an expression of my gratitude.


\section{Stratified spaces}\label{secstrat}


This section contains basic material concerning stratifications of manifolds: 
regularity conditions, Whitney stratifications, and normally nonsingular inclusions.
We also recall the notion of a Bekka stratification~\cite{Bekka91} 
and describe the opposite stratifications
to which our main theorem applies.


\subsection{Stratifications}

\begin{defn}
A {\em stratification} of a topological space
$X$ consists of a locally finite collection $\CS=\{S_\alpha\}$ of locally closed subspaces
$S_\alpha\subset X$ called {\em strata} satisfying

\begin{enumerate}
\item
(covering) $X=\bigcup_\alpha S_\alpha$, 
\item
(pairwise disjoint) $S_\alpha \cap S_\beta =\emptyset,$ for $\alpha\not =\beta$, 
\item
(axiom of frontier) 
$\ol {S^c_\alpha} \cap S^c_\beta \not =\emptyset$ if and only if $S^c_\beta \subset \ol {S^c_\alpha}$
where $S_\alpha^c, S_\beta^c$ denote connected components of $S_\alpha,S_\beta$
respectively.
\end{enumerate}
\end{defn}

We would like the group of stratification-preserving homeomorphisms of 
a space
to act transitively on each of the connected components of its strata.
In order for this to be true, we must impose some regularity on the stratification.
Furthermore, we would like whatever regularity we impose to behave well with respect
to common constructions such as forming products.
There are many approaches to this problem which balance strength and applicability.
In this paper, we follow the development in which differential topology plays the lead role.
In particular, we assume that our topological spaces are manifolds and the strata
are submanifolds. Furthermore, in order to keep technical issues at a minimum,
we work in the subanalytic context. Thus all manifolds are assumed
to be real analytic and all stratifications are assumed to be subanalytic.

\begin{defn} 
Let $\CS=\{S_\alpha\}$, $\CT=\{T_\beta\}$ be stratifications of manifolds $M$, $N$ respectively.
A map $f:M\to N$ is said to be {\em stratified} with respect to $\CS$ and $\CT$ provided that
\begin{enumerate}

\item
for each stratum $T_\beta$, the inverse-image $f^{-1}(T_\beta)$
is a union of connected components $S^c_\alpha$ of strata $S_\alpha$,

\item
if $f$ maps a connected component $S^c_\alpha$ of a stratum $S_\alpha$ to the stratum $T_\beta$,
then the restriction $f|{S^c_\alpha}:S^c_\alpha\to T_\beta$ is a submersion.
\end{enumerate}

\end{defn}

In this paper, whether strata are connected are not will not play any role.
For clarity, from hereon, we shall ignore this issue and
the reader may assume strata are connected or take connected components where necessary.



\subsection{Whitney stratifications} 
We recall here the usual notion of a Whitney stratification.

Let $M$ be a manifold, let $X, Y\subset M$ be 
disjoint submanifolds, and let $x\in X\cap \ol Y$. 
Choose local coordinates around $x$ so that some neighborhood in $M$
is identified with an open subspace of affine space. 
The following conditions are easily seen to be independent of the choice of local coordinates.

\begin{defn} ($A$-regularity)
The triple $(Y, X,x)$ is said to satisfy {\em Whitney's condition $A$} 
if 
given any sequence of points
$y_i\in Y$ converging to $x$, such that the tangent planes $T_{y_i} Y$
converge to some plane $\tau$, we have
$$
T_x X\subset \tau.
$$
We say that the pair $(Y,X)$ satisfies {\em Whitney's condition $A$} if 
for all $x\in X$, the
triples $(Y,X,x)$ satisfy
the condition.
\end{defn}

\begin{defn}[{$B$-regularity}]
The triple $(Y, X,x)$ is said to satisfy {\em Whitney's condition $B$} 
if 
given any sequences of points $x_i\in X$ and 
$y_i\in Y$ each converging to $x$, such that the secant lines $\ell_i=\ol{x_i y_i}$
converge to some line $\ell$ and the tangent planes $T_{y_i} Y$
converge to some plane $\tau$, we have
$$
\ell\subset \tau.
$$
We say that the pair $(Y,X)$ satisfies {\em Whitney's condition $B$} if  for all $x\in X$,
the
triples $(Y,X,x)$ satisfy
the condition.
\end{defn}

\begin{rmk} 
It is easy to check that condition $B$ implies condition $A$.
We have separated the two conditions for later use.
\end{rmk}

\begin{defn}
A stratification $\CS=\{S_\alpha\}$ of a manifold $M$ is said to be a
{\em Whitney stratification} (or {\em $B$-regular})
provided that
each pair of strata $(S_\alpha,S_\beta)$ satisfies Whitney's condition $B$.
\end{defn}

The main impact of a stratification $\CS=\{S_\alpha\}$ of a manifold $M$
being a Whitney stratification is that it then
admits the structure of {\em abstract stratified space} 
in the sense of Thom-Mather~\cite{Thom69,Math70,GWduPL76}.
This consists of a collection of {\em tubes} $\CT_\CS=\{(U_\alpha,\pi_\alpha,\rho_\alpha)\}$,
where $U_\alpha\subset M$ is an open neighborhood of $S_\alpha$,
$\pi_\alpha:U_\alpha\to S_\alpha$ is a retraction (so $\pi_\alpha|S_\alpha=\Id_{S_\alpha}$), and
$\rho_\alpha:U_\alpha\to \BR_{\geq 0}$ is a distance function
(so $\rho_\alpha^{-1}(0)=S_\alpha$.)
The tubes must satisfy
delicate compatibility axioms including that they are {\em good} in the following sense.
 
\begin{defn}
A tube $(U_\alpha,\pi_\alpha,\rho_\alpha)$ is said to be {\em good} if the restriction
$$
(\pi_\alpha,\rho_\alpha):U_\alpha\cap S_\beta\to S_\alpha\times \BR
$$
is a submersion whenever $S_\alpha\subset \ol S_\beta\setminus S_\beta$.  
\end{defn}

Vector fields which are {\em controlled} by the tubes
integrate to give continuous stratification-preserving
motions of the manifold.
For example, to see that a Whitney stratification is equisingular along a stratum,
one constructs a controlled vector field
such that the resulting motion moves along the stratum.
To find such a vector field, one shows that the tubes may be chosen
compatibly with certain stratified maps, and controlled vector fields
may be lifted to controlled vector fields. One way to summarize the technique is in
the two ``isotopy lemmas" of Thom recalled in the next section.


\subsection{Bekka stratifications}

For technical reasons, we will also need the notion of a Bekka stratification which we recall here.
All of the definitions and results of this section and the next
may be found in~\cite{Bekka91}.

Let $f:M\to N$ be a map of manifolds, let $X,Y\subset M$ be 
disjoint smooth submanifolds, and let $x\in X\cap \ol Y$. 
Choose local coordinates around $x$ so that some neighborhood in $M$
is smoothly identified with an open subspace of affine space. 
The following condition is easily seen to be independent of the choice of local coordinates.

\begin{defn}[{$A_f$-regularity}]
The triple $(Y, X,x)$ is said to satisfy {\em Thom's condition $A_f$} 
if 
given any sequence of points
$y_i\in Y$ converging to $x$, such that the planes $\ker(df|T_{y_i} Y)$
converge to some plane $\tau$, we have
$$
\ker(df|T_x X)\subset \tau.
$$
We say that the pair $(Y,X)$ satisfies {\em Thom's condition $A_f$} if for all $x\in X$,
the
triples $(Y,X,x)$ satisfy
the condition.
\end{defn}

\begin{defn} 
Let $\CS=\{S_\alpha\}$, $\CT=\{T_\beta\}$ be stratifications of manifolds $M$, $N$ respectively.
A map $f:M\to N$ is said to be a {\em Thom map}
if 
\begin{enumerate}
\item $f$ is stratified, 
\item for each pair of strata $(S_\alpha,S_\beta)$, $f$ satisfies Thom's condition $A_f$.
\end{enumerate}
\end{defn}


\begin{defn} 
A stratification $\CS=\{S_\alpha\}$ of a manifold $M$ is said to be a
{\em Bekka stratification} (or {\em $C$-regular})
provided that
for each stratum $S_\alpha$, there is an open neighborhood $U_\alpha\subset M$ of $S_\alpha$,
and a map
$\rho_\alpha:U_\alpha\to \BR_{\geq 0}
$
such that 
\begin{enumerate}
\item
$\rho_\alpha^{-1}(0)=S_\alpha$,
\item $\rho_\alpha$ is a Thom map with respect
to the restriction of $\CS$ to $U_\alpha$,
and the stratification of $\BR_{\geq 0}$ by $\{0\}$ and $\BR_{>0}$.
\end{enumerate}
\end{defn}

\begin{prop} 
Any Whitney stratification $\CS$ of a manifold $M$ is a Bekka stratification.
Furthermore, all pairs of strata of a Bekka stratification satisfy Whitney's condition $A$.
To summarize, we have
$$
\mbox{
$B$-regularity $\implies$ $C$-regularity $\implies$ $A$-regularity.
}
$$
\end{prop}

One can make the difference between Whitney and Bekka stratifications precise
by considering what are the good tubes $(U_\alpha,\pi_\alpha,\rho_\alpha)$
and changing a single quantifier. 
Whitney stratifications are characterized by: 
for {\em any} ``tubular function" $\rho_\alpha$ (that is, $\rho_\alpha$ is a sum
of the squares of the normal coordinates in some local chart), and any local retraction $\pi_\alpha$,
there is a neighborhood $U_\alpha$ which forms a good tube~\cite{Trot83}.
Bekka stratifications are characterized by: there {\em exists} a distance
function $\rho_\alpha$ such that for any local retraction $\pi_\alpha$, 
there is a neighborhood $U_\alpha$ which forms a good tube~\cite{Bekka91}.

Although Bekka stratifications are weaker than Whitney stratifications,
they nevertheless allow one to use the Thom-Mather stratification theory. 
In particular, they admit the Thom-Mather structure of abstract stratified space,
and we have the following Thom isotopy lemmas.
Let $P$ be a manifold stratified by the single stratum $P$ itself.

\begin{thm}[{Thom's first isotopy lemma}] \label{t1} Let $M$ be a manifold, $\CS=\{S_\alpha\}$ a Bekka stratification of $M$, 
and $f:M\to P$ a proper stratified map.
Then for any $p\in P$, there is an open neighborhood $U\subset P$ of $p$,
and a (not necessarily differentiable)
stratification-preserving homeomorphism $h$ 
fitting into a commutative diagram
$$
\begin{array}{ccccc}
M & & \stackrel{h}{\longrightarrow} & & f^{-1}(p)\times U \\
& \hspace{-1em}f\searrow \hspace{1em}& &  \hspace{1em}\swarrow \pi\hspace{-1em}& \\
& & U & & 
\end{array}
$$
\end{thm}

\begin{thm}[Thom's second isotopy lemma] \label{t2} Let $M$, $N$
be manifolds, $\CS=\{S_\alpha\}$, $\CT=\{T_\beta\}$ Bekka stratifications of $M$ $N$ respectively,
and $f:M\to N$ a proper Thom map. Let $g:N\to P$ be 
a proper stratified map.
Then for any
$p\in P$, there is an open neighborhood $U\subset P$ of $p$,
and (not necessarily differentiable)
stratification-preserving homeomorphisms $h_M$ and $h_N$ 
fitting into a commutative diagram
$$
\begin{array}{ccccc}
M & & \stackrel{h_M}{\longrightarrow} & & (g\circ f)^{-1}(p)\times U \\
\downarrow & & & & \downarrow \\
N & & \stackrel{h_N}{\longrightarrow} & & g^{-1}(p)\times U \\
& \hspace{-1em}f\searrow \hspace{1em}& &  \hspace{1em}\swarrow \pi\hspace{-1em}& \\
& & U & & 
\end{array}
$$
\end{thm}


\subsection{Transversality} 

Another important property Bekka stratifications share with Whitney stratifications 
is that they behave well with respect
to common constructions.


\begin{defn}
Let $\CS=\{S_\alpha\}$, $\CT=\{T_\beta\}$
be stratifications of manifolds $M$, $N$ respectively.
A map $f:M\to N$ is said to take $\CS$ {\em transversely} to $\CT$ if it takes each
$S_\alpha$ transversely to each $T_\beta$.
\end{defn}

\begin{prop}\label{trans}
Let $M$, $N$
be manifolds, $\CS=\{S_\alpha\}$, $\CT=\{T_\beta\}$ Bekka stratifications of $M$, $N$ respectively,
and $f:M\to N$ a map that takes $\CS$ transversely to $\CT$.
Then 
$$\CS\cap f^{-1}(\CT)=\{S_\alpha\cap f^{-1}(T_\beta)\}$$
is a Bekka stratification of $M$.
\end{prop}



\subsection{Normally nonsingular inclusions}

Let $M$ be a manifold, and let $N\subset M$ be a submanifold.
Let $\CS$ be a Bekka stratification of $M$ transverse to $N$,
and let $N\cap \CS$ be the induced Bekka stratification of $N$.
We would like a stratified description of a small neighborhood of $N$ in $M$.

Let $p:E\to N$ be the normal bundle of $N$ in $M$. 
Given an 
inner product 
$g:\Sym^2 E\to \BR$
on the fibers of $E$, and a small $\epsilon > 0$,
let $E_\epsilon\subset E$ be the open neighborhood of the zero section
consisting of points $e\in E$ with $g(e,e)<\epsilon$.
A {\em tubular neighborhood} of $N$ in $M$ consists of an 
inner product $g$
on the fibers of $E$, a choice of small $\epsilon>0$, and an 
embedding $\phi_\epsilon:E_\epsilon\to M$
which is the identity on $N$. 
If we were unconcerned with the stratification $\CS$, 
then a tubular neighborhood would give a satisfactory description
of a small neighborhood of $N$ in $M$.

%

In the stratified setting, we have the following result.

\begin{thm} \label{nni}
There is an inner product on $E$, a small $\epsilon>0$, 
an embedding $\phi_\epsilon:E_\epsilon\to M$ whose image we denote by $M_\epsilon$,
and a (not necessarily differentiable) homeomorphism
$j_\epsilon:E_\epsilon \to M_\epsilon$ which is the identity on $N$, such that
the induced projection
$$
p_\epsilon :M_\epsilon\stackrel{j_\epsilon^{-1}}{\to} E_\epsilon
\stackrel{p}{\to} N
$$
is stratification-preserving. 
\end{thm}

\begin{proof}
The proof for Whitney stratifications given in~\cite{GMSMT88} works verbatim for Bekka stratifications. 
It uses only: (1) maps $f:M\to N$ of manifolds which are transverse to a fixed Bekka stratification
of $N$ are open in the space of all maps (this is equivalent
to the stratification being $A$-regular~\cite{trot78}),
(2) Thom's first isotopy lemma~\ref{t1}
and (3)~\propref{trans}. 
\end{proof}



\subsection{Opposite stratifications}

We describe here the geometric set-up to which our main construction applies.
There are many possible levels of generalization depending on what kinds of spaces, groups, 
and actions one might consider. In this paper, we restrict ourselves to the specific context
of complex geometry.

Let $M$ be a complex manifold and let $a:\BC^\times \times M\to M$ be an action.
For a point $p\in M$ fixed by $a$, we have the ascending and descending spaces
$$
S^+_p=\{q\in M| \lim_{z\to \infty} a(z)\cdot q = p\}
\qquad
S^-_p=\{q\in M| \lim_{z\to 0} a(z)\cdot q = p\}
$$
By definition, the ascending and descending spaces are contractible.
We write $\CS^+=\{S^+_p\}$ and 
$\CS^-=\{S^-_p\}$ for the collections of ascending and descending spaces respectively.

\begin{defn}
 An action $a:\BC^\times \times M\to M$ is said to be {\em Morse-Whitney} provided that

\begin{enumerate}
\item there are finitely many $\BC^\times$-fixed points $p\in M$, 

\item each fixed point $p\in M$ is contained in a neighborhood $U\subset M$
which is $\BC^\times$-equivariantly isomorphic to a vector space $V$ 
on which $\BC^\times$ acts linearly,

\item $\CS^+$ and 
$\CS^-$ form transverse Whitney stratifications.

\end{enumerate}
Similarly, it is said to be {\em Morse-Bekka} if we only require Bekka stratifications in condition (3).
\end{defn}

\begin{defn}
A (not necessarily regular) stratification $\CS=\{S_\alpha\}$ of $M$
is said to be {\em simple} if
for each $S_\alpha$, 
the frontier $\ol S_\alpha\setminus S_\alpha$
is a Cartier divisor in $S_\alpha$ (or in other words, locally given by a single equation) or empty.

\end{defn}


\begin{example}
In Section~\ref{secflag}, we recall why the action of a regular one-parameter
subgroup $\BC^\times\subset G$ of a complex reductive group on its flag variety $\CB$
is a Morse-Whitney action whose ascending and descending manifold stratifications are simple.
\end{example}


\section{Sheaves}\label{secsh}

In this section, we briefly recall the definitions of perverse sheaves, intersection
cohomology, and tilting sheaves on a stratified space. 
We state some of their basic properties
but refer the reader to~\cite{GMIHII,BBD82,KS94, BBM04} for any discussion.


\subsection{Constructible sheaves}

Let $M$ be a manifold, and let $\BC_M$ denote the sheaf of locally constant
complex-valued functions. Let $D(M)$ be the bounded derived category
of sheaves of $\BC_M$-modules with $\BR$-constructible cohomology. 
Given an object $\CF$ of $D(M)$, we write $H^*(\CF)$ for its cohomology sheaves. 

Let $\CS=\{S_\alpha\}$ be a Bekka stratification of $M$.
We say that an object $\CF$ of $D(M)$ is $\CS$-{constructible} 
if 
for each stratum $S_\alpha$, the restriction $H^*(\CF)|_{S_\alpha}$ is locally constant.
We write $D_\CS(M)$ for the full subcategory of $D(M)$ of $\CS$-constructible objects.

Thanks to the Thom Isotopy Lemmas~\ref{t1} and~\ref{t2} and \propref{trans},
working with constructible objects on Bekka stratifications is no different from
working with them on Whitney stratifications. 
In particular, standard operations preserve constructibility.
A notional convention: all functors in this paper are assumed to be derived,
though our notation will not always reflect this. For example, for a map $f$,
we write $f_*$ in place of $Rf_*$ to denote the derived pushforward.


\subsection{Perverse sheaves}

We
recall here the notions of perverse sheaves, intersection cohomology,
and tilting sheaves. 
Throughout what follows, 
let $\CS=\{S_\alpha\}$ be a
complex Bekka stratification of a complex manifold $M$.
To simplify the discussion, we assume the strata are connected.


\begin{defn} The category $P_\CS(M)$ of {\em perverse sheaves} is the full subcategory of $D_\CS(M)$
consisting of objects $\CF$ such that for each stratum $s_\alpha:S_\alpha\to M$, we have

\begin{enumerate}

\item $H^k(s_\alpha^*\CF)=0$, for all $k>-\dim_\BC S_\alpha$, 

\item $H^k(s_\alpha^!\CF)=0$, for all $k<-\dim_\BC S_\alpha$.
\end{enumerate}

\end{defn}

The category $P_\CS(M)$ is an Artinian abelian category. 

For each stratum $s_\alpha:S_\alpha\to M$, and a local system $\CL_\alpha$ on $S_\alpha$,
we have the {\em standard object}
$\CL_{\alpha*}=s_{\alpha*}\CL_{\alpha}[\dim_\BC S_\alpha]$ and {\em constandard object}
$\CL_{\alpha!}=s_{\alpha!}\CL_{\alpha}[\dim_\BC S_\alpha]$ of the category $D_\CS(M)$.
The following lemma is a well known implication of the Morse theory
of Stein manifolds.

\begin{lem}\label{stein}
If the stratification $\CS=\{S_\alpha\}$ is simple,
the standard object $\CL_{\alpha*}$ and costandard object $\CL_{\alpha!}$ are perverse.
\end{lem}


\subsection{Intersection cohomology}
Although we will not work with  intersection cohomology sheaves, we recall their defining vanishing conditions for comparison with those of tilting sheaves.

\begin{thm} Given a local system $\CL$ on a stratum $s_0:S_0\to M$,
there exists a unique object $\IC_0(\CL)$ of $P_\CS(M)$ such that
\begin{enumerate}
\item[(0)] $s_0^*\IC_0(\CL)\simeq\CL[\dim_\BC S_0]$

\end{enumerate}
 
\noindent
and for each stratum $s_\alpha:S_\alpha\to M$, with $\alpha\not = 0$, we have

\begin{enumerate}

\item[(1)] $H^k(s_\alpha^*\IC_0(\CL))=0$, for all $k\geq-\dim_\BC S_\alpha$, 

\item[(2)] $H^k(s_\alpha^!\IC_0(\CL))=0$, for all $k\leq -\dim_\BC S_\alpha$, 

\end{enumerate}

We call $\IC_0(\CL)$
the {\em intersection cohomology sheaf} of $S_0$ with coefficients in $\CL$.
\end{thm}

An object of $P_\CS(M)$ is simple if and only if it is the intersection
cohomology sheaf of a stratum
with coefficients in an irreducible local system.



\subsection{Tilting sheaves}

\begin{defn} 
An object $\CF$ of $P_\CS(M)$ is a {\em tilting sheaf}
if for each stratum $s_\alpha:S_\alpha\to M$, we have
%
%
\begin{enumerate}

\item $H^k(s_\alpha^*\CF)=0$, for all $k\not =-\dim_\BC S_\alpha$, 

\item $H^k(s_\alpha^!\CF)=0$, for all $k\not =-\dim_\BC S_\alpha$.

\end{enumerate}

\end{defn}

The following results may be found in~\cite{BBM04}.  
The second is a special case
of very general categorical results on tilting objects.
Recall that if $\CS=\{S_\alpha\}$ is a simple stratification, 
then for a local system $\CL_\alpha$ on a stratum $S_\alpha$, 
the standard object  $\CL_{\alpha*}$ and costandard object $\CL_{\alpha!}$ are perverse.

\begin{prop} \label{tiltfilt}
Suppose $\CS=\{S_\alpha\}$ is a simple stratification. Then a perverse sheaf $\CF$
is tilting if and only if both of the following hold
\begin{enumerate}
\item
$\CF$ may be represented as a successive extension
of standard objects
 $\CL_{\alpha*}$,
 \item $\CF$ may be represented as a successive extension
of costandard objects
$\CL_{\alpha!}$.
 \end{enumerate}
 
\end{prop}

We say that a tilting sheaf is {\em minimal} if it is indecomposable in the full subcategory
 of all tilting sheaves.

\begin{prop} Suppose $\CS=\{S_\alpha\}$ is a simple stratification,
and for each $S_\alpha$, we have $\pi_1(S_\alpha)=\pi_2(S_\alpha)=\{1\}$.
Then there is a bijection
between minimal
tilting sheaves and strata. The minimal tilting sheaf $T_\alpha$ corresponding to the stratum 
$S_\alpha$ satisfies
\begin{enumerate}
\item $s_\alpha^*T_\alpha\simeq\BC_{S_\alpha}[\dim_\BC S_\alpha]$,
\item the closure $\ol S_\alpha$ is the support of $T_\alpha$.
\end{enumerate}

\end{prop}

Under the assumptions of the proposition, we have a natural bijection
between strata, simple perverse sheaves, and minimal tilting sheaves.
It is also true that each simple perverse sheaf has a minimal projective cover.


\section{Morse theory}\label{secmorse}

This section contains an account of ideas from~\cite{KS94, GMVMPS}.


\subsection{Critical points}
Let $M$ be a complex manifold. 
The contangent bundle $\pi:T^*M\to M$ is naturally a symplectic manifold.

Given a submanifold $S\subset M$, the {\em conormal bundle} $\pi:T^*_S M\to S$ is
defined by the natural short exact sequence
$$
0\to T^*_S M \to T^*M|_S\to T^*S\to 0.
$$
It is a Lagrangian submanifold.

Given a complex Whitney stratification $\CS=\{S_\alpha\}$, the {\em conormal Lagrangian}
$\Lambda_\CS\subset T^*M$
is defined to be the union 
$$
\Lambda_\CS=\bigcup_\alpha T^*_{S_\alpha}M.
$$
It is a $\BC^\times$-conical Lagrangian subvariety.

\begin{defn} 
A point $p\in M$ is a {\em critical point} of a function $f:M\to \BC$ with respect to the stratification $\CS$
if the differential
$df:M\to T^*M$ takes $p$ to a point of $\Lambda_\CS$.
\end{defn}

The following is immediate from the definitions. 

\begin{lem}
Let $S_\alpha$ be the stratum containing $p$. 
The following are equivalent:
\begin{enumerate}
\item $p$ is not a critical point of $f$,
\item
$df|T_p S_\alpha$ is a submersion,
\item
$df|T_pM$ is nonzero,
and $T_p f^{-1}(f(p))$ is transverse to $T_p S_\alpha$,
\item
$f$ takes $S_\alpha$ at $p$ transversely to the point $f(p)$.
\end{enumerate}
\end{lem}


\subsection{Boundary-stratified balls}

\begin{defn} Let $M$ be a manifold. A {\em boundary-stratified ball} in $M$
is a closed embedded topological ball $B\subset M$ 
equipped with a stratification $\CT$ of its boundary $\partial B$
such that the complement $M\setminus B$, the interior $B\setminus \partial B$, and $\CT$
form a Whitney stratification of $M$.
\end{defn}

It is perhaps misleading that we use the term ``ball," since as the following example
demonstrates, we have in mind other shapes such as cylinders. This
is why we allow the boundary to be stratified rather than assuming it is a submanifold.

\begin{example} \label{boundstratex}
Let $M$ be a manifold, and let $N\subset M$ be a submanifold.
Choose a tube $(U,\pi,\rho)$  around $N$ in $M$. Choose a closed 
embedded ball $B_N \subset N$
such that $\partial B_N\subset N$ is a submanifold. Define $B_\epsilon\subset U$ to be the
inverse image
$$
B_\epsilon=(\pi\times\rho)^{-1}(B_N\times [0,\epsilon]).
$$
Then for sufficiently small $\epsilon> 0$, $B_\epsilon \subset M$ is naturally a boundary-stratified ball
(in fact a manifold with corners)
with boundary stratification $\CT_\epsilon$ consisting of three strata: 
$$
T^v_\epsilon = (\pi\times\rho)^{-1}(\partial B_N\times [0,\epsilon)) 
$$
$$
T^h_\epsilon = (\pi\times\rho)^{-1}((B_N\setminus\partial B_N)\times \{\epsilon\}) 
$$
$$
T^c_\epsilon = (\pi\times\rho)^{-1}(\partial B_N\times \{\epsilon\}) 
$$
The superscript letters stand for ``vertical," ``horizontal," and ``corner."
\end{example}

%




\subsection{Test quadruples} 


Let $M$ be a complex manifold equipped with a fixed complex
Whitney stratification $\CS=\{S_\alpha\}$.

\begin{defn} Consider quadruples of data $Q=(B, f, D, e)$ where

\begin{enumerate}
\item
$B\subset M$ is a boundary-stratified ball with boundary stratification $\CT$, 

\item $f:B\to\BC$ is the restriction of  an analytic function defined on
an open neighborhood of $B$,

\item $D\subset f(B)\subset \BC$ is a smoothly embedded closed disk,

\item $e\in D$ is a point.

\end{enumerate}

\noindent
The quadruple $Q=(B, f, D, e)$ is said to be a {\em test quadruple} if
there is an open neighborhood $V\subset\BC$ containing $D$ such that the data satisfies

\begin{enumerate}

\item in the open set $f^{-1}(V)$, 
the boundary stratification $\CT$ is transverse to the stratification $\CS$,

\item
in the open set $f^{-1}(V)$,
for each stratum $S_\alpha\in \CS$ and each boundary stratum $T_\beta\in \CT$,
the restriction $f|S_\alpha\cap T_\beta$ is a submersion,

\item
the set $\partial D\cup e$ is disjoint from the critical values of $f$.
\end{enumerate}

\end{defn}


\subsection{Test functions} We show here that test quadruples are abundant.

Let $M$ be a complex manifold equipped with a fixed
complex Bekka stratification $\CS=\{S_\alpha\}$.
Fix a point $p\in M$, and let $F$ be the germ of an analytic function at $p$.
So there is a small open neighborhood $U_p\subset M$ containing $p$
such that we have a well-defined function $F:U_p\to \BC$. 
By translation, we may arrange so that $F(p)=0$. 

Suppose that $F$ has an isolated
critical point at $p$. In other words, the 
graph of the differential $dF:U_p\to T^*M$ intersects the conormal Lagrangian 
$\Lambda_\CS\subset T^*M$
in a single covector $\xi\in\Lambda_\CS$, and $\xi$ lies in the fiber of $\Lambda_\CS$ above $p$.

\begin{prop}
There exists a test quadruple $Q_{p,F}=(B,f,D,e)$ such that $p\in B\subset U_p$, and
$f$ is the restriction of $F$.
\end{prop}

\begin{proof}
To begin, recall the construction of a boundary-stratified ball of Example~\ref{boundstratex}.
Choose a good tube $(U_\alpha,\pi_\alpha,\rho_\alpha)$  around the stratum 
$S_\alpha\subset M$ containing $p$.
It is technically convenient to assume that $\pi_\alpha$ and $\rho_\alpha$ are analytic
near $p$ which we may since $\CS$ is a Whitney stratification.
Choose a small closed Euclidean ball $B_\alpha\subset S_\alpha$ containing $p$
with spherical boundary $\partial B_\alpha\subset S_\alpha$.
We have a boundary-stratified ball
$$
B_\epsilon=(\pi_\alpha\times\rho_\alpha)^{-1}(B_\alpha\times [0,\epsilon]),
$$
for any small $\epsilon> 0$, with boundary stratification $\CT_\epsilon$. 
The following confirms requirement (1)
of a test quadruple
for any boundary-stratified ball $B_\epsilon$ of this form.

\begin{lem} The boundary stratification $\CT_\epsilon$
of $B_\epsilon$ is transverse to $\CS$.
\end{lem}

\begin{proof}
Follows immediately from the definition of a good tube
and the functoriality of Whitney stratifications.
\end{proof}

Next, we turn to requirement (2) of a test quadruple.
Since $p$ is an isolated critical point of $F$,
to confirm requirement (2) for some open neighborhood $V\subset \BC$ containing $F(p)=0$,
it suffices to show that for some choice of $B_\alpha$ and $\epsilon>0$, 
the level set $V_0=F^{-1}(0)$ cuts through the strata of $\CS\cap \CT_\epsilon$ 
transversely. 

We first check it for intersections of the strata of $\CT_\epsilon$ 
with the stratum $S_\alpha$ itself. 
 But the only stratum of $\CT_\epsilon$ to meet $S_\alpha$
 is the ``vertical" stratum
 $$
T^v_\epsilon = (\pi_\alpha\times\rho_\alpha)^{-1}(\partial B_\alpha\times [0,\epsilon)),
$$
and we have
$$
\partial B_\alpha = T^v_\epsilon\cap S_\alpha.
$$

\begin{lem} For small enough $B_\alpha\subset S_\alpha$ containing $p$, 
the boundary submanifold
$\partial B_\alpha$ is transverse to $V_0$ inside $M$.
\end{lem}

\begin{proof}
Since $p$ is an isolated critical point, $V_0$ is transverse to $S_\alpha\setminus \{p\}$
inside $M$. And for small enough $B_\alpha$ containing $p$,
it is a standard result that $V_0\cap S_\alpha$ is transverse to $\partial B_\alpha$ inside $S_\alpha$.
\end{proof}

We assume in what follows that we have fixed a suitably small ball $B_\alpha$
so that the previous lemma holds.
Next consider any stratum $S_\beta$ and its intersection
with  
the ``corner" stratum
 $$
T^c_\epsilon = (\pi_\alpha\times\rho_\alpha)^{-1}(\partial B_\alpha\times \{\epsilon\}).
$$

\begin{lem}
For any small enough $\epsilon>0$, $V_0$ cuts through the intersection $S_\beta\cap T^c_\epsilon$
transversely.
\end{lem}

\begin{proof}
By functoriality,
the inverse image $\pi_\alpha^{-1}(\partial B_\alpha) \subset U_\alpha$ is naturally Whitney stratified
so that the restriction $\rho_\alpha|\pi_\alpha^{-1}(\partial B_\alpha)$ 
is a Thom map. 
Together with the previous lemma, this implies the assertion for any small $\epsilon>0$.
\end{proof}

Note that by the lemma, for small enough $\epsilon>0$,
we immediately have that $V_0$ cuts through the intersection $S_\beta\cap T^v_\epsilon$ transversely as well.

Finally, consider the intersection of $S_\beta$
with the ``horizontal" stratum
$$
T^h_\epsilon = 
(\pi_\alpha\times\rho_\alpha)^{-1}((B_\alpha\setminus\partial B_\alpha)\times \{\epsilon\}).
$$
Suppose for all $\epsilon>0$, $V_0$ does not cut 
through $S_\beta\cap T_\epsilon^h$ transversely. Then we arrive at a contradiction
by applying the curve selection lemma.

Finally, we may choose $D\subset \BC$ to be any suffciently small disk around $0$,
and $e$ to be any point in $D$ distinct from $0$.
\end{proof}
%



\subsection{Morse groups} 

Let $X$ be a locally compact topological space, 
$j:Y\to X$ the inclusion of 
a closed subspace,
and $i:U=X\setminus Y\to X$ the inclusion of the open complement.
For any object $\CF$ of $D(X)$, 
we have a natural distinguished triangle
$$
j_! j^! \CF\to \CF\to  i_* i^* \CF\stackrel{[1]}{\to}.
$$
Taking hypercohomology, 
we obtain the long exact sequence 
$$
H^*(X,U;\CF)\to H^*(X;\CF)\to H^*(U;\CF)\stackrel{[1]}{\to}.
$$
For example, if $\CF$ is the dualizing sheaf $\BD_X$, we obtain the usual long exact sequence
in (Borel-Moore) homology 
$$
H_{-*}(Y)\to H_{-*}(X)\to H_{-*}(X,Y)\stackrel{[1]}{\to}.
$$

We adopt the perspective of Morse theory 
that the term $H^*(U;\CF)$ represents
the ``difference" between the other two terms.
In other words, if one starts with $H^*(X,U;\CF)$ and tries to build 
$H^*(X;\CF)$, the amount necessary to be ``added" and ``subtracted"
is $H^*(U;\CF)$, with whether to add or subtract being determined
by the maps.
There are situations when the maps are particularly simple: most notably,
when the long exact sequence is non-zero in only one degree. In this case,
we have a short exact sequence, there is no cancellation, and the middle term
may be realized (non-canonically) as a sum of the two outer terms.
Such miraculous vanishing is at the heart of what it means for
$\CF$ to be a perverse sheaf. 

Let $M$ be a complex manifold equipped with a complex
Whitney stratification $\CS$.
For a test quadruple $Q=(B, f, D, e)$, let 
$$
i_Q: B\cap f^{-1}(D\setminus e) \to M
$$ 
denote the inclusion.

\begin{defn}
To an object $\CF$ of $D_\CS(M)$,
and a test quadruple $Q$, we define the {\em Morse group} $M^*_Q(\CF)$
to be
$$
M^*_Q(\CF) =R\Gamma(M,i_{Q*}i_Q^! \CF)
$$
\end{defn}

The following is a fundamental connection between perverse sheaves and Morse theory.

\begin{thm}
For an object $\CF$ of $D_\CS(M)$, the following are equivalent:
\begin{enumerate}
\item $\CF$ is a perverse sheaf,
\item
for every test quadruple $Q$, we have
$$
M^k_Q(\CF)
= 0, \mbox{ for } k\not = 0,
$$
\item
for every test quadruple $Q=(B,f,D,e)$ such that $df$ takes $B$ transversely to $\Lambda_\CS$, 
we have
$$
M^k_Q(\CF)
= 0, \mbox{ for } k\not = 0.
$$
\end{enumerate}
\end{thm}

\begin{proof}
The equivalence of $(1)$ and $(3)$ is Theorem~10.3.12 and Proposition~6.6.1 of \cite{KS94}. A simple perturbation argument shows $(3)$ implies $(2)$.
\end{proof}


\subsection{Morse kernels}
We introduce here the sheaf which represents the Morse group 
of a test quadruple. 
Recall that for a test quadruple $Q=(B, f, D, e)$, we write
$$
i_Q: B\cap f^{-1}(D\setminus e) \to M
$$ 
for the inclusion.

\begin{defn}
For a test quadruple $Q$, the {\em Morse kernel} $\CM_Q$
is the object of $D(M)$ defined by
$$
\CM_Q=i_{Q!} i_Q^* \BC_M
$$
\end{defn}

\begin{prop}
The Morse kernel $\CM_Q$ represents the Morse group $M^*_Q$.
To be precise,
for any object $\CF$ of $D_\CS(M)$, we have 
$$
M^*_Q(\CF) = 
 R\Hom_{D(M)}^*(\CM_Q,\CF).
$$
\end{prop}

\begin{proof}
By adjunction, we have
$$
R\Gamma(M,i_{Q*}i_Q^! \CF)\simeq R\Hom_{D(M)}(\BC_M,i_{Q*} i_Q^!\CF)
\simeq R\Hom^*_{D(M)}(i_{Q!} i_Q^*\BC_M,\CF).
$$
\end{proof}

We will need the following for our main construction.

\begin{prop} 
For a test quadruple $Q=(B, f, D, e)$, the Morse kernel $\CM_Q$
is constructible with respect to a Whitney stratification 
transverse to $\CS$.
\end{prop}

\begin{proof}
%
%
%
%
%
Choose an open subspace $V\subset \BC$ containing $D$ such that
the requirements of the definition of a test quadruple are satisfied.
Let $\CT$ be the boundary stratification of $B$.
Then $V$ admits a Whitney stratification $\CV$ with the following stratum closures:
the 
point $e$, the boundary of the disk $\partial D$, the disk $D$,
and $V$ itself.
The restriction of the boundary stratification $\CT$ 
provides a Whitney stratification $\CU$ of the inverse image $U=f^{-1}(V)$.
It follows from the requirements of a test quadruple that the map $f$ takes $\CU$
transversely to $\CV$. Therefore by transversality, the stratification 
$$
\CW=f^{-1}(\CV)\cap \CU
$$
is a Whitney stratification. Since the construction of the Morse kernel $\CM_Q$
involves only $\CW$-constructible objects and standard functors,
it is $\CW$-constructible.

It remains to show that $\CW$ is transverse to $\CS$.
Again by the requirements of a test quadruple,
it follows that $f$ takes $\CS$
transversely to $\CV$.
Therefore the decomposition 
$
f^{-1}(\CV)\cap \CS
$
is a Whitney stratification. It suffices to show that the boundary stratification $\CT$
is transverse
to $f^{-1}(\CV)\cap \CS$. We consider the problem in two regions.
First, in $f^{-1}(V \setminus (\partial D\cup e))$, we immediately
have the desired assertion by the requirements of a test quadruple
and the fact that in this region the strata of $f^{-1}(\CV)$ are open.
Second, using requirement (3), 
choose a small neighborhood $V'\subset V$ containing $\partial D\cup e$
such that $f$ does not have any critical values in $V'$.
By requirement (2), in the region $f^{-1}(V')$, $f$ takes $\CT\cap \CS$ transversely to $\CV$.
This is a simple reformulation of the desired assertion.
\end{proof}


\section{From Morse theory to tilting sheaves}\label{secmain}


\subsection{Nearby cycles}

Let $M$ be a complex manifold
equipped with a 
(not necessarily complex) Bekka stratification $\CS$. 

Let $D\subset \BC$ be the open unit disk stratified
by the origin $0\in D$ and its complement $D^\times\subset D$.

Let $\pi:M\to D$ be a proper nonsingular map.
Let $M^\times\subset M$ denote the subspace $\pi^{-1}(D^\times)$,
and $M_0\subset M$ the fiber $\pi^{-1}(0)$.

We assume that $\pi$ is {stratified} with repect to the given stratifications.
Recall that by definition, this means
\begin{enumerate}
\item $M_0$ and $M^\times$ are each a union of strata, which we denote by $\CS_0$
and $\CS^\times$,

\item $\pi$ maps each stratum of $\CS^\times$ submersively to $D^\times$.
\end{enumerate}

With this set-up, we have the nearby cycles functor
$$
R\psi_{\pi}:D_{\CS^\times}(M^\times)\to D_{\CS_0}(M_0).
$$

\begin{caution} If we assume the stratification $\CS$ is complex,
then it is well known that the nearby cycles functor restricts to a functor on
perverse sheaves
$$
R\psi_{\pi}:P_{\CS^\times}(M^\times)\to P_{\CS_0}(M_0)
$$
But in what follows, we shall be considering stratifications $\CS$ which are not
complex and taking nearby cycles of objects which are not perverse. Reinforcing 
this potential confusion, we will nevertheless obtain perverse sheaves in this way.
\end{caution}


\subsection{Main construction}

Let $M$ be a compact complex manifold equipped with a $\BC^\times$-action 
$$
a: \BC^\times \times M\rightarrow M.
$$
Consider the diagonal $\BC^\times$-action on $M\times \BC$, where $\BC^\times$ acts on $M$
via the given action $a$, and on $\BC$ via the standard multiplication. 

Suppose $M\times \BC$ is equipped with a $\BC^\times$-invariant Bekka stratification $\CS$. In other words, each stratum of $\CS$ is preserved by the $\BC^\times$-action.
Then the natural projection 
$$
\pi:M\times \BC\to\BC
$$ 
is a stratified map
with respect to $\CS$ and the stratification of $\BC$ by $0$ and $\BC^\times$.
Thus we have the nearby cycles functor
$$
R\psi_{\pi}:D_{\CS^\times}(M\times\BC^\times)\to D_{\CS_0}(M)
$$
where $\CS^\times$ denotes the restriction of $\CS$ to the subspace $M^\times=M\times \BC^\times$, and $\CS_0$ the restriction to the fiber $M_0=M$.

Now let $\CS_1$ denote the restriction of $\CS$ to the fiber $M_1=\pi^{-1}(1)$ which we also 
identify with $M$.
Consider the projection
$$
\tilde a: M\times \BC^\times\rightarrow M
$$
given by the inverse of the action
$$
\tilde a(m,z)=a(z^{-1})\cdot m.
$$
By construction, the pullback functor preserves constructibility
$$
\tilde a^*:D_{\CS_1}(M)\to D_{\CS^\times}(M).
$$

\begin{defn}
Given an object $\CF$ of $D_{\CS_1}(M)$, we define the
object $\Psi(\CF)$ of $D_{\CS_0}(M)$ to be
the nearby
cycles
$$
\Psi(\CF)= R\psi_{\pi} \tilde a^*\CF
$$
where $\pi:M\times\BC\to\BC$ is the natural projection,
and $\tilde a:M\times\BC^\times\to M$ is the inverse of the action.
\end{defn}

\begin{rmk} \label{invtriv}
If $\CS_1$ is $\BC^\times$-invariant, then
for any object $\CF$ of $D_{\CS_1}(M)$, we
have a natural isomorphism
$$
\Psi(\CF)\simeq\CF.
$$
\end{rmk}


\subsection{Constructing stratifications}
In the preceding discussion, we assumed we were given 
a $\BC^\times$-invariant Bekka stratification $\CS$
of the product $M\times \BC$ equipped with the product $\BC^\times$-action.
Here we describe an interesting way in which such stratifications arise.

Recall that a $\BC^\times$-action on $M$is said to be Morse-Bekka if the fixed points are isolated,
the ascending manifolds $\CS^+$ and descending
manifolds $\CS^-$ form transverse Bekka stratifications of $M$,
and each fixed point has a neighborhood on which the action is linear.

\begin{defn}
Given a Morse-Bekka $\BC^\times$-action on $M$, and a
Bekka stratification $\CS_1$ of $M$
such that $\CS_1$
is {transverse} to the descending manifold stratification $\CS^-$. 
We define the {\em flow stratification} $\CS$ of the family $M\times\BC$
as follows. In the subspace $M^\times=M\times\BC^\times$, we take $\CS$
to be the result $\CS^\times$ of acting on $\CS_1$ via the diagonal $\BC^\times$-action.
In the fiber $M_0=M$, we take $\CS$ to be the ascending manifold stratification $\CS^+$.
\end{defn}

\begin{prop}  The flow stratification $\CS$ is a $\BC^\times$-invariant 
Bekka stratification of the family $M\times\BC$.
\end{prop}

\begin{proof}
By construction, $\CS$ is $\BC^\times$-invariant, and its restriction to the subspace $M^\times$ 
is a Bekka stratification.
We must check that $\CS$ is $C$-regular at the strata $\CS^+$ of the zero fiber $M_0$.
The proof is by induction beginning with open strata and following the 
closure containment partial order. 
The initial case and each of the succesive inductive steps is implied by the following lemma.
\end{proof}

\begin{lem} Let $V$ be a vector space on which $\BC^\times$ acts linearly
without any trivial summands. In other words,
$V$ is the direct sum $V_+\oplus V_-$, where $V_+\subset V$ denotes the subspace
on which the weights of the action are positive, and 
$V_-\subset V$ the subspace on which the weights are negative.
Let $\CS$ be a $\BC^\times$-invariant Bekka stratification of $V\setminus V_+$
which is transverse to $V_-$. Then $\CS$ together with $V_+$ forms a Bekka stratification
of $V$.
\end{lem}

\begin{proof}
Follows immediately from more general results of Bekka~\cite{Bekka97}.
\end{proof}

\begin{rmk}
We do not know if the lemma is true for Whitney stratifications. This is our primary reason
for working in the generality of Bekka stratifications.
\end{rmk}

\begin{defn}
Given a Morse-Bekka $\BC^\times$-action on $M$,
we say that an object $\CF$ of $D(M)$ is in {\em good position} if it is 
constructible with respect to a Bekka stratification $\CT$ transverse 
to the decending manifold stratification $\CS^-$.
\end{defn}
\begin{cor}\label{constr}
If an object $\CF$ of $D(M)$ is in good position
with respect to a Morse-Bekka $\BC^\times$-action,
then the object $\Psi(\CF)$ is constructible with respect to the ascending manifold stratification $\CS^+$.
\end{cor}


\subsection{Constant kernels} 

In this section, we return to the viewpoint that an object $\CK$ of $D(M)$
provides a contravariant functor
$
R\Hom^*_{D(M)}(\CK,\cdot).
$

\begin{prop}\label{constker}
If an object $\CK$ of $D(M)$ is 
is in good position
with respect to a Morse-Bekka $\BC^\times$-action,
then for any object $\CF$ of $D_{\CS^-}(M)$, we have a natural isomorphism
$$
R\Hom^*_{D(M)}(\Psi(\CK),\CF)\simeq R\Hom^*_{D(M)}({\CK},\CF).
$$
\end{prop}

\begin{proof}
Let $\CT$ be a Bekka stratification of $M$ transverse to the descending manifold
stratification $\CS^-$ such that $\CK$ is $\CT$-constructible.

Consider the fiber product 
$$
\pi^{(2)}:M\times M\times\BC\to \BC
$$ 
of two copies of the family
$\pi:M\times\BC\to\BC$.
We have the fiber product Bekka stratification 
$$
\CS^{(2)}=\CS\underset{\BC}{\times}\CS^-
$$
where $\CS$ is the flow decomposition associated to $\CT$,
and $\CS^-$ is the {constant} extension of the descending manifold stratification.
By construction, the stratifications $\CS$ and $\CS^-$ are transverse,
so the relative diagonal
$$
\Delta_{M\times \BC}\subset M\times M\times \BC
$$
is transverse to the fiber product stratification $\CS^{(2)}$.
Therefore $\Delta_{M\times \BC}$ is Bekka stratified by its intersection with
$\CS^{(2)}$. Furthermore, by Theorem~\ref{nni}, there is a small neighborhood 
$$
\CN_{M\times\BC}\subset M\times M\times \BC
$$
of $\Delta_{M\times \BC}$ together with a retraction
$$
r:\CN_{M\times\BC}\to\Delta_{M\times\BC}
$$
such that with respect to
the stratification of $\CN_{M\times\BC}$ induced by $\CS^{(2)}$, the retraction $r$ is a stratified 
fiber bundle
whose fibers are open balls. 

We conclude that for an object $\CF^{(2)}$ of $D_{\CS^{(2)}}(M\times M\times\BC^\times),$
we have a natural isomorphism
$$
 i^* R\psi_{\pi^{(2)}}(\CF^{(2)})
\simeq
R\psi_{\pi} i^* (\CF^{(2)})
$$
where we write $i:\Delta_{M\times \BC}\to M\times M\times \BC 
$
for the inclusion. 

Now we calculate. In what follows, to make the formulas more readable, 
we write $p:M\to pt$ for the augmentation map.
Rewriting $R\Hom^*$ in terms of $\otimes$ and $\verdier$, we have
$$
R\Hom^*({\Psi(\CK)},\CF) 
= 
p_*R\inthom(\Psi(\CK), \CF)
\simeq 
p_*\verdier(\verdier(\CF)\otimes\Psi(\CK)).
$$
By definition and Remark~\ref{invtriv}, we have
$$
p_*\verdier(\verdier(\CF)\otimes\Psi(\CK))
\simeq
p_*\verdier(\verdier(R\psi_\pi \tilde a^*\CF)\otimes R\psi_\pi \tilde a^*\CK).
$$
By a standard identity for $\otimes$, we have
$$
p_*\verdier(\verdier(R\psi_\pi \tilde a^*\CF)\otimes R\psi_\pi \tilde a^*\CK)
\simeq
p_*\verdier (i^*(\verdier (R\psi_\pi \tilde a^*\CF)\boxtimes R\psi_\pi \tilde a^*\CK)).
$$
Now we apply standard identities and the isomorphism deduced above to obtain
$$
p_*\verdier (i^*(\verdier (R\psi_\pi \tilde a^*\CF)\boxtimes R\psi_\pi \tilde a^*\CK))
\simeq
R\psi_{z} \pi_*\verdier (i^*(\verdier (\tilde a^*\CF)\boxtimes \tilde a^*\CK)).
$$
Here we have written $z$ for the standard coordinate on $\BC$,
and $R\psi_z$ for the associated nearby cycles functor.
Rewriting in terms of $R\inthom$, we have
$$
R\psi_{z} \pi_*\verdier (i^*(\verdier (\tilde a^*\CF)\boxtimes \tilde a^*\CK))
\simeq
R\psi_{z} \pi_* R\inthom( \CK,\CF)\boxtimes\BC_{\BC^\times}.
$$
Finally, unwinding the definitions, we have 
$$
R\psi_{z} \pi_* R\inthom( \CK,\CF)
\simeq 
R\psi_{z} (R\Hom^*( \CK,\CF)\boxtimes\BC_{\BC^\times})
\simeq 
R\Hom^*(\CF,\CK).
$$
\end{proof}



\subsection{Constructing tilting sheaves} 

Starting with an object $\CK$ of $D(M)$ that is in good position
with respect to a Morse-Bekka $\BC^\times$-action,
we would like to analyze the stalks and costalks of the object $\Psi(\CK)$.
By Corollary~\ref{constr}, we know that $\Psi(\CK)$ is constructible with respect to
the ascending manifold stratification $\CS^+=\{S^+_\alpha\}$.
Therefore, given such an object,  it suffices to calculate its stalks and costalks
at the fixed points $\{p_\alpha\}$ of the $\BC^\times$-action.

For each $\alpha$, let $i_\alpha:p_\alpha\to M$ be the inclusion of the fixed point, and
let $s^+_\alpha:S^+_\alpha\to M$ and 
$s^-_\alpha:S^-_\alpha\to M$ be the inclusions of the ascending and descending manifolds
respectively.

\begin{lem}
For any object $\CK$ of $D_{\CS^+}(M)$, we have canonical isomorphisms
$$
i_{p_\alpha}^* s^{+*}_\alpha\CK\simeq R\Gamma(M,s_{\alpha*}^-s^{-*}_\alpha\CK)
$$
$$
i_{p_\alpha}^* s^{+!}_\alpha\CK\simeq R\Gamma(M,s_{\alpha!}^-s^{-*}_\alpha\CK)
$$
\end{lem}

\begin{proof} Since
$S^-_\alpha$ is a contractible
transverse slice 
to $S^+_\alpha$,
both assertions are standard implications of Theorem~\ref{nni}.
\end{proof}


Now recall the standard and costandard objects of $D_{\CS^-}(M)$ defined by
$$
\CJ_{\alpha*}=s^-_{\alpha*}\BC_{S^-_\alpha}[\dim_\BC S_\alpha^-]
$$
$$
\CJ_{\alpha!}=s^-_{\alpha!}\BC_{S^-_\alpha}[\dim_\BC S_\alpha^-]
$$

\begin{lem}
For any object $\CK$ of $D_{\CS^+}(M)$, we have canonical isomorphisms
$$
i_{p_\alpha}^* s_\alpha^{+*}\CK\simeq  \verdier(R\Hom_{D(M)}^*(\CK,\CJ_{\alpha!}))
[-\dim_\BC S_\alpha^-]
$$
$$
i_{p_\alpha}^* s_\alpha^{+!}\CK\simeq  \verdier(R\Hom_{D(M)}^*(\CK,\CJ_{\alpha*}))
[-\dim_\BC S_\alpha^-]
$$
\end{lem}

\begin{proof}
For the first assertion, rewriting $R\Hom^*$ in terms of $\otimes$ and $\verdier$, we have
$$
R\Hom_{D(M)}^*(\CK,\CJ_{\alpha!}) \simeq 
R\Gamma(M,\verdier(\verdier(s^-_{\alpha!}\BC_{S^-_\alpha})\otimes\CK))[\dim_\BC S^-_\alpha] 
$$
and then by standard identities involving $\verdier$ and pushforwards, we have 
$$
R\Gamma(M,\verdier(\verdier(s^-_{\alpha!}\BC_{S^-_\alpha})\otimes\CK))[\dim_\BC S^-_\alpha] 
\simeq 
\verdier(R\Gamma_c(M,s^-_{\alpha*}\BD_{S^-_\alpha}\otimes\CK))[\dim_\BC S^-_\alpha].
$$
Since $S_\alpha^-$ is smooth and canonically oriented (since it is complex), we have
$$
\verdier(R\Gamma_c(M,s^-_{\alpha*}\BD_{S^-_\alpha}\otimes\CK))[\dim_\BC S^-_\alpha] 
\simeq 
\verdier(R\Gamma_c(M,s^-_{\alpha*}\BC_{S^-_\alpha}\otimes\CK))[-\dim_\BC S^-_\alpha] 
$$
We have the identity
$$
\verdier(R\Gamma_c(M,s^-_{\alpha*}\BC_{S^-_\alpha}\otimes\CK))[-\dim_\BC S^-_\alpha] 
\simeq 
\verdier(R\Gamma_c(M,s^-_{\alpha*}s^{-*}_\alpha\CK))[-\dim_\BC S^-_\alpha] 
$$
Since $M$ is compact, we are done by the previous lemma. 

The second assertion is similar. 
\end{proof}

We may now state our main result. 
Recall that a $\BC^\times$-action on $M$ is said to be Morse-Whitney if the fixed points are isolated,
the ascending manifolds $\CS^+$ and descending
manifolds $\CS^-$ form transverse Whitney stratifications of $M$,
and each fixed point has a neighborhood on which the action is linear.

\begin{thm} \label{main}
Let $M$ be a manifold with a Morse-Whitney $\BC^\times$-action.
Suppose that the descending manifold stratification $\CS^-$ is a simple stratification.
If $\CM$ is a Morse kernel with respect to $\CS^-$, 
then $\Psi(\CM[\dim_\BC M])$ is a tilting perverse sheaf 
on the ascending manifold stratification $\CS^+$.
\end{thm}

\begin{proof}
We have seen that $\Psi(\CM[\dim_\BC M])$ is $\CS^+$-constructible
so it remains to see its stalks and costalks are in the correct degree.
For this, recall that by Proposition~\ref{constker}, we have an isomorphism 
$$
R\Hom_{D(M)}^*(\Psi(\CM[\dim_\BC M]),\CF)\simeq 
R\Hom_{D(M)}^*(\CM[\dim_\BC M]),\CF)
$$ for any $\CS^-$-constructible object $\CF$ . Thus if $\CF$ is perverse,
we have 
$$
R\Hom_{D(M)}^*(\Psi(\CM[\dim_\BC M]),\CF)\simeq V[-\dim_\BC M]
$$
where $V$ is the Morse group of $\CF$ with respect to $\CM$ placed in degree $0$. 
Applying this to the standard object $\CJ^-_{\alpha !}$
and costandard object $\CJ^-_{\alpha *}$,
the previous lemma implies
$$
i_{p_\alpha}^* s_\alpha^{+*}\Psi(\CM[\dim_\BC M])\simeq \verdier(V_!)[\dim_\BC S^+_\alpha]
$$
$$
i_{p_\alpha}^* s_\alpha^{+!}\Psi(\CM[\dim_\BC M])\simeq \verdier(V_*)[\dim_\BC S^+_\alpha]
$$
where $V_!$ and $V_*$ are the Morse groups of $\CJ^-_{\alpha !}$
and $\CJ^-_{\alpha *}$. Since the stratification $\CS^-$ is simple,
$\CJ^-_{\alpha !}$
and $\CJ^-_{\alpha *}$ are perverse, and so $V_!$ and $V_*$ are concentrated
in degree $0$. Thus the stalk and costalk
cohomology at $p_\alpha$ is in the correct degree, and since 
$\Psi(\CM[\dim_\BC M])$ is $\CS^+$-constructible, we are done.
\end{proof}


%





\section{Flag varieties}\label{secflag}

Let $G$ be a connected reductive complex algebraic group.
Let $B^+,B^-\subset G$ be opposite Borel subgroups, and let 
$\frakb^+,\frakb^-$ be the Lie algebras of $B^+, B^-$ respectively.
Let $N^+, N^-\subset B^+$ be the unipotent radicals of $B^+, B^-$ respectively, and
let $\frakn^+,\frakn^-$ be the Lie algebras of $N^+, N^-$ respectively.
Let $T\subset G$
be the maximal torus $T=B^+\cap B^-$.
Let $W=N_G(T)/T$ be the Weyl group of $G$, and 
let $w_0\in W$ be the longest element.

Let $\CB$ be the flag variety of $G$. 
Choose a coweight $\lambda:\BC^\times\to T$ so that the weights of $\BC^\times$
acting on $\frakn^-$ via $\lambda$
are positive (and therefore the weights of $\frakn^+$ are negative).
Let 
$
a:\BC^\times\times \CB\to\CB
$ 
be the action
defined by composing $\lambda$ with the standard $G$-action.
We recall briefly some standard results which confirm that $a$ is a good action
in the sense of Section~\ref{secstrat}.

The fixed points of $a$ are naturally indexed by $w\in W$:
if $\tilde w\in N_G(T)$ represents $w$, then the fixed point $p_w\in\CB$
labeled by $w$ is given by $\tilde w\cdot B^+$.
So in particular, we have
$p_1=B^+$ and $p_{w_0}=B^-$. 
For each fixed point $p_w\in \CB$,
the ascending and descending manifolds 
$$S^+_w=\{p\in\CB| \lim_{z\to \infty} a(z)\cdot p = p_w\}
\qquad
S^-_w=\{p\in\CB| \lim_{z\to 0} a(z)\cdot p = p_w\}$$
coincide with the $N^+$ and $N^-$-orbits through $p_w$ respectively.
Thus the ascending and descending manifold
stratifications $\CS^+=\{S^+_w\}$ and $\CS^-=\{S^-_w\}$ are Whitney stratifications.
Since $S^+_w$ and $S^-_w$ also coincide with the $B^+$ and $B^-$-orbits through $p_w$
respectively, and $\frakb^+$ and $\frakb^-$ are transverse subalgebras,
the stratifications $\CS^+$ and $\CS^-$ are transverse.

For each fixed point $p_w\in \CB$, consider the Lie algebra $\frakn_{ww_0}$ 
of the unipotent radical $N_{ww_0}$ of the opposite Borel subgroup $B_{ww_0}\subset G$.
Let $\CU_w\in\CB$ be the open neighborhood of $p_w$ consisting of all Borel subgroups
opposite to $B_{ww_0}$.
The action of $\frakn_{ww_0}$ on $p_w$ via the exponential map and the standard $G$-action
provides an isomorphism $\frakn_{ww_0}\risom \CU_w$. 
Moreover, the isomorphism interwines the linear $\BC^\times$-action on $\frakn_{ww_0}$ via $\lambda$ and
the $\BC^\times$-action on $\CU_w$ via $a$.
Thus we have a locally linear model at each fixed point.

Finally, for each $w\in W$, the inclusions of strata $s_w^+:S_w^+\to \CB$ and 
$s_w^-:S_w^-\to \CB$
are affine, and so the stratifications $\CS^+$ and $\CS^-$ respectively are simple.


\subsection{Equivariantization}

Let $D_{N^+}(\CB)$, $D_{N^-}(\CB)$ denote the
$N^+$,$N^-$-equivariant bounded derived categories of $\CB$ respectively.
Since $N^+$ and $N^-$ are contractible,
the forgetful functors provide natural equivalences 
$$
F^+:D_{N^+}(\CB)\risom D_{\CS^+}(\CB)
\qquad
F^-:D_{N^-}(\CB)\risom D_{\CS^-}(\CB)
$$
In what follows, we use these identifications without further comment.

We also have the equivariantization functors
$$
\Gamma^+:D(\CB)\risom D_{N^+}(\CB)
\qquad
\gamma^-:D(\CB)\risom D_{N^-}(\CB)
$$
defined as follows. Using the diagrams
$$
\CB\stackrel{p_+}{\leftarrow} N^+\times\CB\stackrel{q^+}{\to}\CB
\qquad
\CB\stackrel{p_-}{\leftarrow} N^-\times\CB\stackrel{q^-}{\to}\CB
$$
where $p$ is the obvious projection, and $q$ is the action, we have
$$
\Gamma^+(\CF)=q^+_*p_+^*(\CF)
\qquad
\gamma^-(\CF)=q^-_!p_-^!(\CF)
$$
They satisfy the following~\cite{MUV92}:
\begin{enumerate}
\item $(F^+,\Gamma^+)$ and $(\gamma^-,F^-)$ are pairs of adjoint functors,

\item the restrictions
$$
\Gamma^+:D_{\CS^-}(\CB)\to D_{\CS^+}(\CB)
\qquad
\gamma^-:D_{\CS^+}(\CB)\to D_{\CS^-}(\CB)
$$ 
are inverse equivalences of categories,

\item for $w\in W$, there are natural isomorphisms
$$
\Gamma^+(\CJ^-_{w,!})\simeq\CJ^+_{w, *}
\qquad
\gamma^-(\CJ^+_{w,*})\simeq\CJ^-_{w, !}
$$
intertwining standard and costandard objects.
\end{enumerate}
In the second statement, we have omitted the forgetful functors from the notation.

Recall that an object $\CF$ of $D(\CB)$ is said to be in good position with respect to the
action 
$$a:\BC^\times \times\CB\to\CB$$ 
if it is constructible with respect to a Bekka stratification $\CT$ transverse to
the descending manifold stratification $\CS^-$.

\begin{prop}\label{equiv}
For any object $\CF$ of $D(\CB)$ which is in good position with respect
to the $\BC^\times$-action $a$, we have an isomorphism
$$
\gamma^-(\CF)\simeq\gamma^-(\Psi(\CF)).
$$
\end{prop}

\begin{proof}
By Proposition~\ref{constker},
for any object $\CF$ of $D_{\CS^-}(M)$, we have a natural isomorphism
$$
R\Hom^*_{D(M)}(\Psi(\CK),\CF)\simeq R\Hom^*_{D(M)}({\CK},\CF).
$$
Hence by adjunction, we have a natural isomorphism
$$
R\Hom^*_{D_{S^-}(M)}(\gamma^-(\Psi(\CK)),\CF)\simeq R\Hom^*_{D_{S^-}(M)}(\gamma^-({\CK}),\CF).
$$
Therefore the assertion follows from the Yoneda lemma.
\end{proof}

\subsection{Projective sheaves} It is well-known that any exact functor on $P_{\CS^-}(\CB)$
is representable by a projective object. 
The following explicitly identifies
the representing object for local Morse groups.

\begin{prop}\label{repproj}
For a Morse kernel $\CM_Q$, its equivariantization $\gamma^-(\CM_Q)$
is the projective object of $P_{\CS^-}(\CB)$ which represents the local Morse group $M_Q$.

\end{prop}

\begin{proof}
For any object $\CF$ of $D_{\CS^-}(\CB)$, by adjunction, we have 
$$
R\Hom^*_{D_{\CS^-}(\CB)}(\gamma^-(\CM),\CF)\simeq R\Hom^*_{D(\CB)}(\CM,\CF).
$$
Thus it suffices to show that $\gamma^-(\CM)$ is perverse. 
We give two arguments, the first geometric and the second categorical.

(1) Geometric: by \propref{equiv}, 
it suffices to show that $\gamma^-(\Psi(\CM))$ is perverse. 
By Theorem~\ref{main}, $\Psi(\CM)$ is a tilting sheaf,
so by Proposition~\ref{tiltfilt}, it is filtered by standard objects $\CJ^+_{w,*}$.
Since $\gamma^-$ takes the standard objects $\CJ^+_{w,*}$ 
to the costandard objects 
$\CJ^-_{w,!}$, it takes tilting sheaves to perverse sheaves. 

(2) Categorical: we know that there exists a projective perverse sheaf $P^-_Q$
which represents $M_Q$. Since $P_{\CS^-}(\CB)$ generates  
$D_{\CS^-}(\CB)$,
we have an isomorphism
$$
R\Hom^*_{D_{\CS^-}(\CB)}(\gamma^-(\CM),\CF)\simeq R\Hom^*_{D_{\CS^-}(\CB)}(P^-_Q,\CF).
$$
Thus we are done by the Yoneda lemma.
\end{proof}

\begin{cor}
For $w\in W$, we have
$$
\gamma^-(T^+_w[-\dim_\BC\CB])\simeq P^-_w
\qquad
\Gamma^+(P^-_w)\simeq T^+_w[-\dim_\BC\CB]
$$
\end{cor}

\begin{proof}
It is easy to see that every indecomposable tilting sheaf $T^+_w$ is a summand
of $\Psi(\CM_Q[\dim_\BC\CB])$ for some Morse kernel. 
Thus by Propositions~\ref{equiv} and~\ref{repproj}, we have
that $\gamma^-$ takes tilting sheaves to projectives. Using that $\gamma^-$
and $\Gamma^+$ are inverse equivalences of categories, and proceeding by induction,
it is easy to check that $\gamma^-$ takes the indecomposable $T^+_w[-\dim_\BC\CB]$ to 
the indecomposable $P^-_w$.
\end{proof}


\bibliographystyle{plain}
\bibliography{ref}

\begin{thebibliography}{10}

\bibitem{BBM04}
A.~Beilinson, R.~Bezrukavnikov, and I.~Mirkovi{\'c}.
\newblock Tilting exercises.
\newblock {\em Mosc. Math. J.}, 4(3):547--557, 782, 2004.

\bibitem{BBD82}
A.~A. Be{i}linson, J.~Bernstein, and P.~Deligne.
\newblock Faisceaux pervers.
\newblock In {\em Analysis and topology on singular spaces, I (Luminy, 1981)},
  pages 5--171. Soc. Math. France, Paris, 1982.

\bibitem{Bekka91}
K.~Bekka.
\newblock C-r\'egularit\'e et trivialit\'e topologique.
\newblock In {\em Singularity theory and its applications, Part I (Coventry,
  1988/1989)}, volume 1462 of {\em Lecture Notes in Math.}, pages 42--62.
  Springer, Berlin, 1991.

\bibitem{Bekka97}
K.~Bekka.
\newblock Regular quasi-homogeneous stratifications.
\newblock In {\em Stratifications, singularities and differential equations, II
  (Marseille, 1990; Honolulu, HI, 1990)}, volume~55 of {\em Travaux en Cours},
  pages 1--14. Hermann, Paris, 1997.

\bibitem{ENV04}
M.~Emerton, D.~Nadler, and K.~Vilonen.
\newblock A geometric {J}acquet functor.
\newblock {\em Duke Math. J.}, 125(2):267--278, 2004.

\bibitem{GMVMPS}
S.~Gelfand, R.~MacPherson, and K.~Vilonen.
\newblock Microlocal perverse sheaves.
\newblock {\em Preprint math.AG/0509440}, 2005.

\bibitem{GWduPL76}
Christopher~G. Gibson, Klaus Wirthm{\"u}ller, Andrew~A. du~Plessis, and Eduard
  J.~N. Looijenga.
\newblock {\em Topological stability of smooth mappings}.
\newblock Springer-Verlag, Berlin, 1976.

\bibitem{GMIHII}
Mark Goresky and Robert MacPherson.
\newblock Intersection homology. {I}{I}.
\newblock {\em Invent. Math.}, 72(1):77--129, 1983.

\bibitem{GMSMT88}
Mark Goresky and Robert MacPherson.
\newblock {\em Stratified {M}orse theory}, volume~14 of {\em Ergebnisse der
  Mathematik und ihrer Grenzgebiete (3)}.
\newblock Springer-Verlag, Berlin, 1988.

\bibitem{KS94}
Masaki Kashiwara and Pierre Schapira.
\newblock {\em Sheaves on manifolds}, volume 292 of {\em Grundlehren der
  Mathematischen Wissenschaften}.
\newblock Springer-Verlag, Berlin, 1994.

\bibitem{Math70}
John Mather.
\newblock {\em Notes on topological stability}.
\newblock Harvard University, 1970, available at
  {www.math.princeton.edu/facultypapers/mather/}.

\bibitem{MUV92}
I.~Mirkovi{\'c}, T.~Uzawa, and K.~Vilonen.
\newblock Matsuki correspondence for sheaves.
\newblock {\em Invent. Math.}, 109(2):231--245, 1992.

\bibitem{Ringel91}
Claus~Michael Ringel.
\newblock The category of modules with good filtrations over a quasi-hereditary
  algebra has almost split sequences.
\newblock {\em Math. Z.}, 208(2):209--223, 1991.

\bibitem{SoeCharForm98}
Wolfgang Soergel.
\newblock Character formulas for tilting modules over {K}ac-{M}oody algebras.
\newblock {\em Represent. Theory}, 2:432--448 (electronic), 1998.

\bibitem{Thom69}
R.~Thom.
\newblock Ensembles et morphismes stratifi\'es.
\newblock {\em Bull. Amer. Math. Soc.}, 75:240--284, 1969.

\bibitem{trot78}
D.~J.~A. Trotman.
\newblock Stability of transversality to a stratification implies {W}hitney
  {$(a)$}-regularity.
\newblock {\em Invent. Math.}, 50(3):273--277, 1978/79.

\bibitem{Trot83}
David Trotman.
\newblock Comparing regularity conditions on stratifications.
\newblock In {\em Singularities, Part 2 (Arcata, Calif., 1981)}, volume~40 of
  {\em Proc. Sympos. Pure Math.}, pages 575--586. Amer. Math. Soc., Providence,
  R.I., 1983.

\end{thebibliography}

\end{document}